\input amsppt.sty
\hsize=12.5cm
\vsize=540pt
\baselineskip=16truept
\NoRunningHeads
\TagsOnRight
 \define\CC{\Bbb C}
 \define\CO{\Cal O}
 \define\eps{\varepsilon}
 \define\halfskip{\vskip6pt}
 \font\msbm=msbm10
 \define\MX{\text{\msbm X}}
 \define\NN{\Bbb N}
 \define\pder#1#2{\frac{\partial#1}{\partial#2}}
 \redefine\phi{\varphi}
 \define\proj{\operatorname{pr}}
 \define\PSH{\Cal P\Cal S\Cal H}
 
 \define\Reg{\operatorname{Reg}}
 \define\skipaline{\vskip12pt}
 \define\too{\longrightarrow}
 \define\wdht{\widehat}
 \define\wdtl{\widetilde}
\catcode`@=11
\def\qed{\ifhmode\textqed\else\ifinner\quad\square
   \else\eqno\square\fi\fi}
\def\textqed{{\unskip\nobreak\penalty50
    \quad\hbox{}\nobreak\hfil$\m@th\square$
    \parfillskip=0pt \finalhyphendemerits=0\par}}
\catcode`@=\active

\document
\topmatter
\title
Cross theorem
\endtitle
\author
Marek Jarnicki (Krak\'ow), Peter Pflug (Oldenburg)
\endauthor

\address
Uniwersytet Jagiello\'nski \newline
Instytut Matematyki \newline
30-059 Krak\'ow, Reymonta 4, Poland\newline
{\it E-mail: }{\rm jarnicki{\@}im.uj.edu.pl}
\endaddress

\address
Carl von Ossietzky Universit\"at Oldenburg
\newline Fachbereich Mathematik\newline
Postfach 2503\newline
D-26111 Oldenburg, Germany\newline
{\it E-mail: }{\rm pflug{\@}mathematik.uni-oldenburg.de}
\endaddress

\abstract
Let $D, G\subset\CC$ be domains, let $A\subset D$, $B\subset G$ be locally
regular sets, and let $X:=(D\times B)\cup(A\times G)$. Assume that $A$ is
a Borel set. Let $M$ be a proper analytic subset of an open neighborhood of $X$.
Then there exists a pure $1$-dimensional analytic subset $\wdht M$ of the
envelope of holomorphy $\wdht X$ of $X$ such that any function separately
holomorphic on $X\setminus M$ extends to a holomorphic function on
$\wdht X\setminus\wdht M$. The result generalizes special cases which were
studied in \cite{\"Okt 1998}, \cite{\"Okt  1999a}, and \cite{Sic 2000}.
\endabstract
\endtopmatter

\noindent{\bf 1. Introduction. Main result.}
For domains $D\subset\CC^n$, $G\subset\CC^m$ and non-pluripolar subsets
$A\subset D$, $B\subset G$, define the {\it cross}
$$
X=\MX(D,A;G,B):=(D\times B)\cup(A\times G)\tag *
$$
(notice that $X$ is connected).
Let $U\subset D\times G$ be an open connected neighborhood of $X$ and let
$M$ be an analytic subset of $U$, $M\neq U$. Put
$$
M_z:=\{w\in G\: (z,w)\in M\},\; z\in D,\quad
M^w:=\{z\in D\: (z,w)\in M\},\; w\in G.
$$
We say that a function $f\:X\setminus M\too\CC$ is {\it separately
holomorphic} on $X\setminus M$ ($f\in\CO_s(X\setminus M)$) if:
$$
\forall_{z\in A\: M_z\neq G}\:\; f(z,\cdot)\in\CO(G\setminus M_z),\qquad
\forall_{w\in B\: M^w\neq D}\:\; f(\cdot,w)\in\CO(D\setminus M^w).
$$

For an open set $\Omega\subset\CC^n$ and $A\subset\Omega$ put
$$
h_{A,\Omega}:=\sup\{u\:\; u\in\PSH(\Omega),\;u\leq1 \text{ on }
\Omega,\; u\leq0 \text{ on } A\},
$$
where $\PSH(\Omega)$ denotes the set of all functions plurisubharmonic on
$\Omega$. Define
$$
\omega_{A,\Omega}:=\lim_{k\to+\infty}h^\ast_{A\cap\Omega_k,\Omega_k},
$$
where $(\Omega_k)_{k=1}^\infty$ is a sequence of relatively compact open
sets $\Omega_k\subset\Omega_{k+1}\subset\subset\Omega$ with
$\bigcup_{k=1}^\infty\Omega_k=\Omega$ ($h^\ast$ denotes the upper
semicontinuous regularization of $h$). Observe that the definition is
independent of the approximation sequence $(\Omega_k)_{k=1}^\infty$.

For a cross (*) put
$$
\wdht X:=\{(z,w)\in D\times G\: \omega_{A,D}(z)+\omega_{B,G}(w)<1\}.\tag **
$$

We say that a subset $A\subset\CC^n$ is {\it locally pluriregular}
if $h^\ast_{A\cap\Omega,\Omega}(a)=0$ for any $a\in A$ and for any open
neighborhood $\Omega$ of $a$ (in particular, $A\cap\Omega$ is non-pluripolar).
As always, if $n=1$, then we say that $A$ is locally `regular' instead
of `pluriregular'.

\newpage

The main result of the paper is the following

\proclaim{Theorem 1} Let $D, G\subset\CC$ be domains, let $A\subset D$,
$B\subset G$ be locally regular sets, and let $X:=(D\times B)\cup(A\times G)$.
Assume that $A$ is a Borel set. Let $M$ be a proper analytic subset of an open
connected neighborhood $U$ of $X$. Then there exists a pure $1$-dimensional
analytic subset $\wdht M$ of $\wdht X$ ($\wdht X$ is given by {\rm (**)})
such that for any $f\in\CO_s(X\setminus M)$ there exists exactly one
$\wdht f\in\CO(\wdht X\setminus\wdht M)$ with $\wdht f=f$ on
$X\setminus(M\cup\wdht M)$.

Moreover, if $U=\wdht X$ and $M$ is pure $1$--dimensional, then
the above condition is satisfied by $\wdht M:=M$.
\endproclaim

\proclaim{Remark} \rm Consider the following general problem. Let
$D_j\subset\CC^{n_j}$ be a domain of holomorhy and let $A_j\subset D_j$ be a
locally
pluriregular Borel set, $j=1,\dots,N$. Define the {\it generalized cross}
$$
X:=(A_1\times\dots\times A_{N-1}\times D_N)\cup\dots\cup
(D_1\times A_2\times\dots\times A_N)\subset
\CC^{n_1}\times\dots\times\CC^{n_N}.
$$
Let $U\subset D_1\times\dots\times D_N$ be a connected neighbourhood of $X$
and let $M\subset U$ be a proper analytic set. A function
$f\:X\setminus M\too\CC$ is said to be {\it separately analytic}
($f\in\CO_s(X\setminus M)$) if for any
$(a_1,\dots,a_N)\in A_1\times\dots\times A_N$
and $k\in\{1,\dots,N\}$ the function
$f(a_1,\dots,a_{k-1},\cdot,a_{k+1},\dots,a_N)$ is holomorphic in the domain
$\{z_k\in D_k\:(a_1,\dots,a_{k-1},z_k,a_{k+1},\dots,a_N)\not\in M\}$. Define
$$
\wdht X:=\{(z_1,\dots,z_N)\in D_1\times\dots\times D_N\:\omega_{A_1,D_1}(z_1)+
\dots+\omega_{A_N,D_N}(z_N)<1\}.
$$
{\bf Conjecture:} There exists a pure $1$--codimensional analytic subset
$\wdht M\subset\wdht X$ such that for any $f\in\CO_s(X\setminus M)$ there
exists an $\wdht f\in\CO(\wdht X\setminus\wdht M)$ with $\wdht f=f$ on
$X\setminus(M\cup\wdht M)$. Moreover, $\wdht M=M$ if $U=\wdht X$ and $M$ is
pure $1$--codimensional.
Compare also \cite{\"Okt 1999b} (for $N=2$ and $U=\wdht X$).

\halfskip

Theorem 1 solves the case $N=2$, $n_1=n_2=1$.

J\. Siciak \cite{Sic 2000} solved the following case: $n_1=\dots=n_N=1$,
$D_1=\dots=D_N=\CC$, $M=P^{-1}(0)$, where $P$ is a non-zero polynomial of
$N$ complex variables; the special subcase $N=2$, $P(z,w):=z-w$ had been
previously studied in \cite{\"Okt 1998}, \cite{\"Okt 1999a}.

The case studied in \cite{Sic 2000} is the only known case with
$n_1+\dots+n_N>2$. In the general case, the answer is not known even if
$U=\wdht X$ and $M$ is pure $1$--codimensional.

\endproclaim

\skipaline

\noindent{\bf 2. Auxiliary results.}
The following lemma gathers a few standard results, which will be used
in the sequel.

\proclaim\nofrills{Lemma 2}\;
{\rm (cf\. \cite{Kli 1991}, \cite{Jar-Pfl 2000}, \S\;3.5).}
{\rm (a)} Let $\Omega\subset\CC^n$ be a bounded open set and let
$A\subset\Omega$. Then:

$\bullet$ If $P\subset\CC^n$ is pluripolar, then
$h^\ast_{A\setminus P,\Omega}=h^\ast_{A,\Omega}$.

$\bullet$ $h_{A_k\cap\Omega_k,\Omega_k}^\ast\searrow h_{A,\Omega}^\ast$
(pointwise on $\Omega$) for any sequence of open sets
$\Omega_k\nearrow\Omega$ and any sequence $A_k\nearrow A$.

$\bullet$ $\omega_{A,\Omega}=h^\ast_{A,\Omega}$.

$\bullet$ The following conditions are equivalent:

\hang{for any connected component $S$ of $\Omega$ the set $A\cap S$ is
non-pluripolar;}

\hang{$h^\ast_{A,\Omega}(z)<1$ for any $z\in\Omega$.}

$\bullet$ If $A$ is non-pluripolar, $0<\alpha<1$, and
$\Omega_\alpha:=\{z\in\Omega\: h^\ast_{A,\Omega}(z)<\alpha\}$,
then for any connected component $S$ of $\Omega_\alpha$ the set
$A\cap S$ is non-pluripolar (in particular, $A\cap S\neq\varnothing$).

\halfskip

\noindent{\rm (b)} Let $\Omega\subset\CC^n$ be an open set and let
$A\subset\Omega$. Then:

$\bullet$ $\omega_{A,\Omega}\in\PSH(\Omega)$.

$\bullet$ If $A$ is locally pluriregular, then $\omega_{A,\Omega}(a)=0$
for any $a\in A$.

$\bullet$ If $P\subset\CC^n$ is pluripolar, then
$\omega_{A\setminus P,\Omega}=\omega_{A,\Omega}$.

$\bullet$ If $A$ is locally pluriregular and $P\subset\CC^n$ is pluripolar,
then $A\setminus P$ is locally pluriregular.

\halfskip

\noindent{\rm (c)} Let $X=\MX(D,A;G,B)$ be a cross as in {\rm (*)}. Then:

$\bullet$ If $A$ and $B$ are locally pluriregular, then $X\subset\wdht X$.

$\bullet$ If $D$ and $G$ are domains of holomorphy, then $\wdht X$
is a region of holomorphy.
\endproclaim

\proclaim{Lemma 3} Let $X=\MX(D,A;G,B)$ be a cross as in {\rm (*)}.
If $A$ and $B$ are locally pluriregular, then $\wdht X$ is a domain.
\endproclaim

\demo{Proof}
It suffices to show that
for any approximation sequences $D_k\nearrow D$, $G_k\nearrow G$ of
relatively compact subdomains with $A\cap D_k\neq\varnothing$,
$B\cap G_k\neq\varnothing$, $k\in\NN$, the sets
$$
\wdht X_k:=\{(z,w)\in D_k\times G_k\: h^\ast_{A\cap D_k,D_k}(z)+
h^\ast_{B\cap G_k,G_k}(w)<1\}, \quad k=1,2,\dots,
$$
are connected. Thus, we may additionally assume that $D$ and $G$ are bounded.
Since the cross $X$ is connected and contained in $\wdht  X$, we only need to
prove that for any $(z_0,w_0)\in\wdht X$, each connected component of the
fiber
$$
\wdht X^{w_0}:=\{z\in D\: (z,w_0)\in\wdht X\}=
\{z\in D\: h^\ast_{A,D}(z)<1-h^\ast_{B,G}(w_0)\}
$$
intersects $A$. If $h^\ast_{B,G}(w_0)=0$, then $\wdht X^{w_0}=D$. If
$h^\ast_{B,G}(w_0)>0$, then we apply Lemma 2(a).
\qed\enddemo

\proclaim\nofrills{Theorem 4}\;{\rm (Classical cross theorem, cf\.
\cite{Ngu-Zer 1991}).} Let $X=\MX(D,A;G,B)$ be as in {\rm (*)}. Assume that:

$\bullet$  $D$, $G$ are domains of holomorphy,

$\bullet$  $A$, $B$ are locally pluriregular,

$\bullet$ $A$ is a Borel set.

\noindent Then for any $f\in\CO_s(X)$ there exists exactly one
$\wdht f\in\CO(\wdht X)$ with $\wdht f=f$ on $X$.
\endproclaim

\proclaim\nofrills{Theorem 5}\;{\rm (Dloussky--Grauert--Remmert theorem,
cf\. \cite{Jar-Pfl 2000}, \S\;3.4).}
Let $\Omega\subset\CC^n$ be a domain and let $M$ be an analytic subset of
$\Omega$. Let $\wdht\Omega$ be the envelope of holomorphy of $\Omega$
(univalent or not).
Then there exists a pure $1$--codimensional analytic subset $\wdht M\subset
\wdht\Omega$ such that for any $g\in\CO(\Omega\setminus M)$ there exists
$\wdht g\in\CO(\wdht\Omega\setminus\wdht M)$ with $\wdht g=g$ on
$\Omega\setminus(M\cup\wdht M)$.

If, moreover, $M=\Omega\cap\wdtl M$, where $\wdtl M$ is a pure $1$--codimensional
analytic subset of $\wdht\Omega$, then  the above condition is satisfied by
$\wdht M:=\wdtl M$.
\endproclaim

\proclaim{Lemma 6} Let $D, G\subset\CC$ be domains, let $A\subset D$,
$B\subset G$ be locally regular sets, and let $X:=\MX(D,A;G,B)$.
Let $M$ be a proper analytic subset of an open connected neighborhood $U$ of
$X$. Assume that $A'\subset A$, $B'\subset B$ are such that:

$\bullet$ $A\setminus A'$ and
$B\setminus B'$ are polar (in particular, $A'$, $B'$ are also locally
regular),

$\bullet$ $M_z\neq G$ for any $z\in A'$ and $M^w\neq D$ for any $w\in B'$.

\noindent {\rm (a)} If $f\in\CO_s(X\setminus M)$ and
$f=0$ on $(A'\times B')\setminus M$, then $f=0$ on $X\setminus M$.

\noindent {\rm (b)} If $g\in\CO(U\setminus M)$ and $g=0$ on
$(A'\times B')\setminus M$, then $g=0$ on $U\setminus M$.
\endproclaim

\demo{Proof} (a) Take a point $(a_0,b_0)\in X\setminus M$. We may assume that
$a_0\in A$. Since $A\setminus A'$ is polar, there exists a sequence
$(a_k)_{k=1}^\infty\subset A'$ such that $a_k\too a_0$. The set
$Q:=\bigcup_{k=0}^\infty M_{a_k}$ is at most countable. Consequently, the set
$B'':=B'\setminus Q$ is non-polar. We have $f(a_k,w)=0$, $w\in B''$,
$k=1,2,\dots$. Hence $f(a_0,w)=0$ for any $w\in B''$. Finally,
$f(a_0,w)=0$ on $G\setminus M_{a_0}\ni b_0$.

\halfskip

(b) Take an $a_0\in A'$. Since $M_{a_0}\neq G$, there exists a
$b_0\in B'\setminus M_{a_0}$. Let
$P=\varDelta_{a_0}(r)\times\varDelta_{b_0}(r)\subset U\setminus M$
($\varDelta_{z_0}(r)$ denote the disc with center $z_0$ and radius $r$). Then
$g(\cdot,w)=0$ on $A'\cap\varDelta_{a_0}(r)$ for any $w\in B'\cap\varDelta_{b_0}(r)$.
The set $A'\cap\varDelta_{a_0}(r)$ is non-polar.
Hence $g(\cdot,w)=0$ on $\varDelta_{a_0}(r)$ for any $w\in B'\cap\varDelta_{b_0}(r)$.
By the same argument for the second variable we get $g=0$ on $P$ and,
consequently, on $U$.
\qed\enddemo

\noindent{\bf 3. Proof of the main theorem.}

\noindent{\bf Step 1.} Fix sequences $D_k\nearrow D$, $G_k\nearrow G$
of relatively compact subdomains with $D_k\subset\subset D_{k+1}$,
$A\cap D_k\neq\varnothing$, $G_k\subset\subset G_{k+1}$,
$B\cap G_k\neq\varnothing$, $k\in\NN$.

For any $a\in A$ such that $M_a\neq G$ we perform the following construction:

Fix a $k\in\NN$, $k\geq2$. Let $M_a\cap G_k=\{b_1,\dots,b_N\}$.
Fix domains $G'=G'_{a,k}$, $G''=G''_{a,k}$ such that
$G_{k-1}\subset\subset G''\subset\subset G'\subset\subset G_k$ and
$b_1,\dots,b_N\in G''$. Take positive numbers $\delta$,
$\eps$, $\eta>\eps$
such that

$\varDelta_a(\delta)\subset\subset D$,

$\varDelta_{b_j}(\eta)\subset\subset G''$, $j=1,\dots,N$,

$\overline\varDelta_{b_i}(\eta)\cap\overline\varDelta_{b_j}(\eta)=\varnothing$,
$i,j=1,\dots,N$, $i\neq j$,

$M\cap(\varDelta_a(\delta)\times G')\subset \bigcup_{j=1}^N\varDelta_a(\delta)
\times\varDelta_{b_j}(\eps)$,

$B\cap V''\neq\varnothing$, where $V'':=G''\setminus\bigcup_{j=1}^N
\overline\varDelta_{b_j}(\eta)$.

Define $V':=G'\setminus\bigcup_{j=1}^N\overline\varDelta_{b_j}(\eps)$.
Note that $V''\subset\subset V'$. Consider the cross
$$
Y=Y_{a,k}:=\MX(\varDelta_a(\delta), A\cap\varDelta_a(\delta); V', B\cap V').
$$
Fix an $f\in\CO_s(X\setminus M)$. Then $f\in\CO_s(Y)$.
By Theorem 4, the function $f$ extends holomorphically to
$\wdht Y\supset\{a\}\times V'$. Consequently, there exists
$0<\wdht\delta<\delta$ such that $f$ is holomorphic
in $\varDelta_a(\wdht\delta)\times V''$.

\halfskip

\noindent{\bf Step 2.} Suppose that for some $j\in\{1,\dots,N\}$ we have:
$$
M\cap(\varDelta_a(\delta)\times\varDelta_{b_j}(\eps))\subset
\{(z,\phi_j(z))\: z\in\varDelta_a(\delta)\},
$$
where $\phi_j\:\varDelta_a(\delta)\too\varDelta_{b_j}(\eps)$ is holomorphic.

We will prove that for sufficiently small $\delta'>0$ the function $f$
extends holomorphically to
$(\varDelta_a(\delta')\times\varDelta_{b_j}(\eta))\setminus
\{(z,\phi_j(z))\: z\in\varDelta_a(\delta')\}$.

Indeed, by Step 1, there exists $\eta'>\eta$ such that the function $f$
extends holomorphically to $\varDelta_a(\wdht\delta)\times(\varDelta_{b_j}(\eta')
\setminus\overline\varDelta_{b_j}(\eta))$. Using the biholomorphism
$$
\varDelta_a(\delta)\times\CC\ni (z,w)\too
(z,w-\phi_j(z))\in\varDelta_a(\delta)\times\CC,
$$
we reduce the problem to the case where $\phi_j\equiv0$. Thus we have
the following problem:

Let $\varDelta(r):=\varDelta_0(r)$.
Given a function $f$ holomorphic on $\varDelta_a(\rho)\times P$, where
$P:=\varDelta(R)\setminus\overline\varDelta(r)$, and such that
$f(z,\cdot)\in\CO(\varDelta(R)\setminus\{0\})$ for any
$z\in A\cap\varDelta_a(\rho)$, prove that $f$ extends
holomorphically to $\varDelta_a(\rho)\times(\varDelta(R)\setminus\{0\})$.

For, consider the cross
$$
Y:=\MX(\varDelta_a(\rho), A\cap\varDelta_a(\rho);\varDelta(R)\setminus\{0\},P).
$$
By Theorem 4, the function $f$ extends to $\wdht Y$.
It remains to observe that
$\wdht Y=\varDelta_a(\rho)\times(\varDelta(R)\setminus\{0\})$
(because $h^\ast_{P,\varDelta(R)\setminus\{0\}}\equiv0$).

\halfskip

In particular, if 
$$
M\cap(\varDelta_a(\delta)\times\varDelta_{b_j}(\eps))=
\{(z,\phi_j(z))\: z\in\varDelta_a(\delta)\},
$$
where $\phi_j\:\varDelta_a(\delta)\too\varDelta_{b_j}(\eps)$ is holomorphic,
for all $j=1,\dots,N$,
then there exists $\delta'>0$ such that $f$ extends holomorphically to
$(\varDelta_a(\delta')\times G'')\setminus M
\supset(\varDelta_a(\delta')\times G_{k-1})\setminus M$.

\halfskip

\noindent{\bf Step 3.} Suppose that for some $j\in\{1,\dots,N\}$ we have:
$$
M\cap(\varDelta_a(\delta)\times\varDelta_{b_j}(\eps))=\{(a,b_j)\}.
$$
By Step 2
(with $\phi_j\equiv b_j$) the function $f$ extends holomorphically to
$\varDelta_a(\delta')\times(\varDelta_{b_j}(\eta)\setminus\{b_j\})$ for some
small $\delta'>0$. On the other hand, we know that
$f$ is separately holomorphic on
$Z:=(\varDelta_a(\delta)\setminus\{a\})\times\varDelta_{b_j}(\eps)$.
Consequently, $f$ is holomorphic on $Z$. Hence $f$ is holomorphic on
$(\varDelta_a(\delta')\times\varDelta_{b_j}(\eps))\setminus\{(a,b_j)\}$. Thus
$(a,b_j)$ is a removable singularity of $f$.

\halfskip

In virtue of the above remark, we may assume that $M$ is pure
$1$-dimensional.

\halfskip

\noindent{\bf Step 4.} Let $A'$ denote the set of all $a\in A$
such that for each $k\geq2$ either there exists $\delta>0$ such that
$M\cap(\Delta_a(\delta)\times G_k)=\varnothing$ or the construction from
Step 1 may be performed in such a way that for each $j\in\{1,\dots,N\}$
$$
M\cap(\varDelta_a(\delta)\times\varDelta_{b_j}(\eps))=
\{(z,\phi_j(z))\: z\in\varDelta_a(\delta)\},
$$
where $\phi_j\:\varDelta_a(\delta)\too\varDelta_{b_j}(\eps)$ is holomorphic
(cf\. Step 2). Then $A\setminus A'$ is at most countable.
Indeed, write
$$
M=\bigcup_{j=1}^\infty\{(z,w)\in P_j\: g_j(z,w)=0\},
$$
where $P_j\subset\subset U$ is a polydisc and
$g_j$ is a defining function for $M\cap P_j$ (cf\. \cite{Chi 1989}, \S\; 2.9).
Put $S_j:=\{(z,w)\in P_j\: g_j(z,w)=\pder{g_j}{w}(z,w)=0\}$.
Observe that if $(z_0,w_0)\in (M\cap P_j)\setminus S_j$, then
there exists a small polydisc $Q=Q'\times Q''\subset\subset P_j$ with center
at $(z_0,w_0)$ such that $M\cap Q$ is the graph of a holomorphic function
$\phi\:Q'\too Q''$.

The projection $\proj_z(S_j)$ is at most countable.
Indeed, we only need to prove that $\proj_z(S'_j)$ is at most countable,
where $S'_j$ is the union of $1$--dimensional irreducible components of
$S_j$. Let $S$ be such an irreducible component. We will show that
$S$ projects onto one point. Take $(z_1,w_1), (z_2,w_2)\in S$.
We want to show that $z_1=z_2$. It suffices to consider only the case where
$(z_1,w_1), (z_2,w_2)$ are regular points of $S$. Let
$\psi=(\psi_1,\psi_2)\:[0,1]\too\Reg(S)$ be a $\Cal C^1$--curve
with $\psi(0)=(z_1,w_1)$, $\psi(1)=(z_2,w_2)$.
Note that $\pder{g_j}{z}(z,w)\neq0$ for $(z,w)\in\Reg(S)$ (because $g_j$ is a
defining function). We have:
$$
0=\pder{(g_j\circ\psi)}{t}(t)=\pder{g_j}{z}(\psi(t))\psi'_1(t),
\quad t\in[0,1].
$$
Thus $\psi'_1\equiv0$. In particular, $z_1=z_2$.

Consequently,
$A\setminus A'\subset\bigcup_{j=1}^\infty\proj_z(S_j)$ is at most countable.

\halfskip
\noindent{\bf Step 5.} Let $B'$ be constructed analogously to $A'$
with respect to the second variable. Put $X':=\MX(D,A';G,B')$.

By Step 2 (and Lemma 6), for any $k\in\NN$ and any
$\xi=(a,b)\in(A'\cap D_k)\times(B'\cap G_k)$
there exists $\rho=\rho_{\xi,k}>0$ such that for each
$f\in\CO_s(X\setminus M)$ there exists $\wdtl f=\wdtl f_{\xi,k}
\in\CO(\Omega_{\xi,k}\setminus M)$
with $\wdtl f=f$ on $X\cap\Omega_{\xi,k}\setminus M$, where
$$
\align
\Omega_{\xi,k}:&=\MX(D_k,\varDelta_a(\rho);G_k,\varDelta_b(\rho))\\
&=(\varDelta_a(\rho)\times G_k)\cup(D_k\times\varDelta_b(\rho))
\subset U\cap(D_k\times G_k).
\endalign
$$
We may always assume that $\rho_{\xi,k+1}\leq\rho_{\xi,k}$.
By Lemma 6, $\wdtl f_{\xi,k+1}=\wdtl f_{\xi,k}$ on
$\Omega_{\xi,k+1}\cap\Omega_{\xi,k}\setminus M$. Define
$$
\Omega:=\bigcup_{k=1}^\infty\bigcup_{\xi\in(A'\cap D_k)\times(B'\cap G_k)}
\Omega_{\xi,k}.
$$

It is clear that $\Omega$ is a connected neighborhood of $X'$.
We will show that the functions $\wdtl f_{\xi,k}$,
$\xi\in(A'\cap D_k)\times(B'\cap G_k)$, $k\in\NN$,
can be glued together. We only need to check that
$\wdtl f_{\xi,k}=\wdtl f_{\eta,k}$ on
$\Omega_{\xi,k}\cap\Omega_{\eta,k}\setminus M$, $\xi=(a,b)$, $\eta=(c,d)$.
Let $\rho':=\rho_{\xi,k}$, $\rho'':=\rho_{\eta,k}$,
$f':=\wdtl f_{\xi,k}$, $f'':=\wdtl f_{\eta,k}$. Observe that
$$
\align
\Omega_{\xi,k}\cap\Omega_{\eta,k}&
=\Big(\varDelta_a(\rho')\times\varDelta_d(\rho'')\Big)\\
&\cup\Big(\varDelta_c(\rho'')\times\varDelta_b(\rho')\Big)\\
&\cup\Big((\varDelta_a(\rho')\cap\varDelta_c(\rho''))\times G_k\Big)\\
&\cup\Big(D_k\times(\varDelta_b(\rho')\cap\varDelta_d(\rho''))\Big)\\
&=:W_1\cup W_2\cup W_3\cup W_4.
\endalign
$$

To prove that $f'=f''$ on $W_1\setminus M$ it
suffices to observe that $f'=f''$
on $(A'\cap\varDelta_a(\rho'))\times(B'\cap\varDelta_d(\rho''))\setminus M$
(and use Lemma 6). The same argument solves the problem on $W_2\setminus M$.

If $W_3\neq\varnothing$, then the equality holds on a non-empty set
$W_3\cap W_1\setminus M$ and we only need to use the identity principle.
The same argument works on $W_4\setminus M$.

\halfskip
\noindent{\bf Step 6.}
Recall that the sets $A'$, $B'$ are locally regular and $A'$ is a Borel set.
Moreover, $h^\ast_{A',D}=h^\ast_{A,D}$ and $h^\ast_{B',G}=h^\ast_{B,G}$.
Hence $\wdht X'=\wdht X$.

First we prove that $\wdht X$ is the envelope of holomorphy of $\Omega$.
We only need to show that any function $g\in\CO(\Omega)$ extends
holomorphically to $\wdht X$. Fix a $g\in\CO(\Omega)$. By Theorem 4
(applied to the cross $X'$), there
exists a $\wdht g\in\CO(\wdht X)$ (recall that $\wdht X=\wdht X'$)
such that $\wdht g=g$ on $X'$. By Lemma 6, $\wdht g=g$ on $\Omega$.

By Theorem 5 there exists a pure $1$--dimensional
analytic subset $\wdht M$ of $\wdht X$ such that for any $g\in\CO(\Omega
\setminus M)$ there exists a $\wdht g\in\CO(\wdht X\setminus\wdht M)$
with $\wdht g=g$ on $\Omega\setminus(M\cup\wdht M)$. We also know that
if $U=\wdht X$ and $M$ is pure  $1$--dimensional, then we can take
$\wdht M=M$.

Now take an $f\in\CO_s(X\setminus M)$ and let $\wdtl f\in\CO(\Omega\setminus
M)$ be such that $\wdtl f=f$ on $X'\setminus M$ (Step 5).
Let $\wdht f\in\CO(\wdht X\setminus\wdht M)$ be
such that $\wdht f=\wdtl f$ in $\Omega\setminus(M\cup\wdht M)$.
In particular, $\wdht f=f$ on $X'\setminus(M\cup\wdht M)$.
By Lemma 6, $\wdht f=f$ on $X\setminus(M\cup\wdht M)$.

Using once again Lemma 6, we conclude that the function
$\wdht f$ is uniquely determined.

\Refs

{\redefine\bf{\rm}
\widestnumber{\key}{XXXXXXXX}
\ref
\key Chi 1989
\by E\. M\. Chirka
\book Complex Analytic Sets
\publ Kluwer Acad\. Publishers
\yr 1989
\endref
\ref
\key Jar-Pfl 2000
\by  M\. Jarnicki, P\. Pflug
\book  Extension of Holomorphic Functions
\publ  de Gruyter Expositions in Mathematics 34, Walter de Gruyter
\yr 2000
\endref
\ref
\key Kli 1991
\by M\. Klimek
\book Pluripotential Theory
\publ Oxford University Press
\yr 1991
\endref
\ref
\key Ngu-Zer 1991
\by Nguyen Thanh Van \& A\. Zeriahi
\paper  Une extension du th\'eor\`eme de
Hartogs sur les fonctions s\'epar\'ement analytiques
\inbook Analyse Complexe Multivariables, R\'ecents D\`evelopements,
A\. Meril (ed.), EditEl, Rende
\yr 1991
\pages 183--194
\endref
\ref
\key \"Okt 1998
\by O\. \"Oktem
\paper Extension of separately analytic functions and applications to
range characterization of exponential Radon transform
\jour Ann\. Polon\. Math\.
\vol 70
\yr 1998
\pages 195--213
\endref
\ref
\key \"Okt 1999a
\by O\. \"Oktem
\paper Extension of separately analytic functions and applications to
mathematical tomography
\jour Dep\. Math\. Stockholm Univ\. (Thesis)
\yr 1999
\endref
\ref
\key \"Okt 1999b
\by O\. \"Oktem
\paper Extending separately analytic functions in $\CC^{n+m}$ with
singularities
\jour Preprint
\yr 1999
\endref
\ref \key Sic 2000
\by J\. Siciak
\paper Holomorphic functions with singularities on algebraic sets
\jour Preprint
\yr 2000
\endref
}
\endRefs

\enddocument